\def\obs#1{{\bf (*** #1 ***)} }

 \def\obs#1{}     
\documentclass[twoside,letterpaper,draft,11pt]{amsart}

\usepackage{amsmath,amsthm,latexsym,xspace,amscd,amssymb,xypic, amscd}               
\usepackage{color}
\usepackage{fullpage}

\newtheorem{teo}{Theorem}[section]
\newtheorem{defi}[teo]{Definition}

\newtheorem{cor}[teo]{Corollary}
\newtheorem{prop}[teo]{Proposition}

\newtheorem{exe}[teo]{Example}

\newcommand{\X}{{\mathbb X}}
\newcommand{\Y}{{\mathbb Y}}

\newcommand{\ex}{{\exists}}

\newcommand{\m}{{}^{-1}}
\newcommand{\mt}{\mapsto}

\def\ndv{\ {\mid \kern -0.7 em {\scriptstyle \not}} \ \ }
\def\nd{\ {\mid \kern -0.4 em {\scriptstyle \not}} \ \ }

\newcommand{\N}{{\mathbb N}}

\numberwithin{equation}{section}
\title[Partial actions of Hausdorff Groups on metric spaces ]{On the continuity of partial actions of Hausdorff groups on metric spaces
}

\begin{document}

\author[J.\ G\'omez ]{J. G\'omez }
\address{Escuela de Matem\'aticas, Universidad Industrial de Santander, Cra. 27 Calle 9 UIS Edificio 45\\  Bucaramanga, Colombia}\email{jorge.gomez14@correo.uis.edu.co }

\author[H.\ Pinedo ]{H. Pinedo }
\address{Escuela de Matem\'aticas, Universidad Industrial de Santander, Cra. 27 Calle 9 UIS Edificio 45\\  Bucaramanga, Colombia}\email{hpinedot@uis.edu.co }
\author[C.\ Uzcategui ]{C. Uzc\'ategui }
\address{Escuela de Matem\'aticas, Universidad Industrial de Santander, Cra. 27 Calle 9 UIS Edificio 45\\  Bucaramanga, Colombia}\email{cuzcatea@uis.edu.co}
\date{\today}

\thanks{The second and third authors thank La Vicerrector\'ia de Investigaci\'on y Extensi\'on de la Universidad Industrial de Santander for the financial support for this work,  which is a part  of the VIE project \# 5761.}
\begin{abstract}  We provide a sufficient condition for a topological partial action of a Hausdorff group on a metric space is continuous, provide that it is separately continuous.
\end{abstract}


\maketitle
\section{Introduction}
The notion of partial action of a group is a weakening of the classical concept of group action and was introduced  in \cite{E1} and \cite{MC}, and was then developed in \cite{AB} and \cite{KL}, in which the authors provided  examples in different guises. Since then partial actions have been an  important tool in $C^* $-algebras and dynamical systems, and in the developing of new cohomological theories \cite{DNP}, \cite{DKh}  and \cite{P}. Every partial action of a group $G$ on a set $\X$ can be obtained, roughly speaking, as a restriction  of a global action (see \cite{AB} and \cite{KL}) on some bigger set $\X_G,$ called the enveloping space of $\X$. Nevertheless, in the category of topological spaces,  when $G$ acts partially on a space $\X$, the superspace  $\X_G$ does not necessarily inherit its topological properties; for instance the globalization of a partial action of a group on a Hausdorff space is not in general Hausdorff (see e.g. \cite[Example 1.4]{AB},  \cite[Proposition 1.2]{AB}). 


On the other hand, actions of Polish groups have important  connections with many areas of mathematics (see \cite{KEB}, \cite{GA} and the references therein). Recall that a Polish space is a topological space which is separable and completely metrizable, and a Polish group is a topological group whose topology is Polish. Partial actions of Polish groups have been recently  considered in  the works \cite{GPU1,PU,PU2}.

It is known that an action of a Polish group $G$ on a metric space $\X$ is continuous provided that it is separately continuous (see \cite[Theorem 3.1.4]{GA}). In this note, under a mild restriction, we generalize  \cite[Theorem 3.1.4]{GA} in two directions, first we only assume that $G$ is Hausdorff and Baire, and second we only assume that $G$ acts partially on $\X.$  

\section{The notions}

Let $G$ be a topological  group with identity element $1$ and $\X$ a topological space. A partially defined map  $G\times \X\to \X$ is a map whose domain is a subset of $G\times \X.$ Let 
$m\colon  G\times \X\rightarrow \X,\,\,(g,x)\mt m(g,x)=g\cdot x\in \X$ be a partial map,   that is $m$ is a map whose domain is contained in $G\times \X.$ We write $\ex g\cdot x$ to mean that $(g,x)$ is in the domain of $m.$
 Then, following \cite{AB, KL} $m$  defines  a (set theoretic) {\it partial action} of $G$ on $\X,$ if for all $g,h\in G$ and $x\in \X$ the following assertions hold:
\smallskip

\noindent (PA1) $\ex g\cdot x$ implies $\ex g\m\cdot (g\cdot x)$ and $g\m\cdot (g\cdot x)=x,$
\smallskip

\noindent (PA2)  $\ex g \cdot (h\cdot x)$ implies $\ex (g h)\cdot x$ and  $g \cdot (h\cdot x)=(g h)\cdot x,$
\smallskip

\noindent (PA3) $\ex 1\cdot x,$ and $1\cdot x=x$.
 \smallskip

Let $G*\X =\{(g,x)\in G\times \X\mid \ex g\cdot x \}$ be the domain of $m,$ set $\X_{g\m}=\{ x\in \X\mid \ex g\cdot x\},$ and $m_g\colon \X_{g\m}\ni x\mt g\cdot x\in\X_g.$
Then a partial action $m\colon G*\X\to\X$ induces a family of bijections $\{m_g\colon \X_{g\m}\to \X_g\}_{g\in G},$ and we denote $m=\{m_g\colon \X_{g\m}\to \X_g\}_{g\in G}.$  
\begin{prop} \label{fam}\cite[Lemma 1.2]{QR} A partial action $m$ of $G$ on $\X$ is a family $m=\{m_g\colon \X_{g\m}\to\X_g\}_{g\in G},$ where $\X_g\subseteq \X,$  $m_g\colon\X_{g\m}\to \X_g$ is bijective, for  all $g\in G,$  and such that:
\begin{itemize}
\item[(i)]$\X_1=\X$ and $m_1=\rm{Id}_\X;$
\item[(ii)]  $m_g( \X_{g\m}\cap \X_h)=\X_g\cap \X_{gh};$
\item[(iii)] $m_gm_h\colon \X_{h\m}\cap  \X_{ h\m g\m}\to \X_g\cap \X_{gh},$ and $m_gm_h=m_{gh}$ in $ \X_{h\m}\cap  \X_{h\m g\m};$
\end{itemize}
for all $g,h\in G.$
\end{prop}
Notice that conditions (ii) and (iii) in Lemma \ref{fam} say that $m_{gh}$ is an extension of $m_gm_h$ for all $g,h\in G.$
We consider   $ G\times \X$ with the product topology and
the subset   $G* \X$ of  $G\times \X$ inherits the  subspace topology. 

\begin{defi} A topological partial action of the group $G$ on the topological space $\X$ is a partial action $m=\{m_g\colon \X_{g\m}\to\X_g\}_{g\in G}$ on the underlying
set $\X$, such that each $\X_g$ is open in $\X$, and each $m_g$ is a homeomorphism, for any $g\in G$. If $m\colon G* \X\to \X$ is continuous, we say that the partial action is continuous. 
\end{defi}


\begin{exe}\label{indu} {\bf Induced partial action:} Let $G$ be a topological group, and $\Y$ a topological space and let $u \colon G\times \Y\to \Y$ be a
 continuous action of  $G$ on $\Y$ and $\X\subseteq \Y$  be an open set. For $g\in G,$ set $\X_g=\X\cap u_g(\X)$ and let $m_g=u_g\restriction \X_{g\m}$ ( the restriction of $u_g$ to $\X_{g\m}$).  Then $m\colon G* \X\ni (g,x)\mt m_g(x)\in \X $ is a topological partial action of $G$ on $\X.$
\end{exe}

%

The interested reader may consult others examples of topological partial actions in \cite[Example 1.2, Remark 1.1, Example 1.3, Example 1.4]{AB} and \cite[p. 108]{KL}.

\section{Topological partial actions of  Hausdorff groups on metric  spaces} Let  $m\colon G*\X\to \X$ be a  topological partial action. For $x\in \X,$ we set $G^x=\{g\in G\mid (g,x)\in G*\X\},$ then by (PA3) $1\in G^x,$ we also set  $m^x\colon G^x\ni g\to m(g,x)\in \X,$ $m$ is called {\it separately continuous} if the maps $m^x$ are continuous, for all $x\in \X.$\footnote{Notice that for any $g\in G$  the map $m_g$ is continuous by definition of topological partial action.}
It is known that group actions of Polish groups on metric spaces are continuous, if and only if, they are separately continuous.  We provide a mild condition on $G^x, \, x\in \X,$ for which separately continuous partial actions of Hausdorff-Baire groups are continuous, the proof we present  is inspired by the one given in   \cite[Theorem 3.1.4]{GA}, for classical group actions.  

\begin{teo}\label{sep} 
Let G be a Hausdorff group, $(\X,d)$ a metric space, and $m$ a  topological partial action of $G$ on $\X.$ Suppose that $G$ is Baire $G^x$ is open in $G,$ for any $x\in \X.$ Then $m$ is continuous, if and only if, it is separately continuous. 
\end{teo}

\proof It is clear that continuous partial actions are separately continuous. For the converse,  suppose that $m$ is separately continuous and let $(g_0,x_0)\in G*\X,$ we check that $m$ is continuous at  $(g_0,x_0).$ Let $l,n\in\N$ and  set
$$
F_{n,l}=\{g \in G^{x_0} \mid \,\forall x \in \X_{g\m} (  d(x, x_0) <2^{-n}\Rightarrow d(m(g, x), m(g, x_0)) \leq 2^{-l})\}.
$$ 
We shall check that that $F_{n,l}$ is a closed subset of $G^{x_0}.$ Indeed,
let $g\in G^{x_0}$ and  $\{g_i\}_{i\in I}$ a net with  $g_i\rightarrow g$ such that  $g_k\in F_{n,l},$ for all $k\in \mathbb{N}.$   Then
	$$
	(\forall k\in \mathbb{N})\left(\forall x\in \X_{g_i^{-1}}\right)\left(d(x,x_0)<\frac{1}{2^n}\Rightarrow
	d\left(m(g_i,x),m(g_i,x_0)\right)\leq\frac{1}{2^l}\right).
	$$
	Let $x\in \X_{g\m},$  that is $g\in G^x,$ since   $G^x$ is open one may assume that $\{g_i\}_{i\in I}\subseteq G^x,$  and $x\in \X_{g\m}\cap \X_{g\m_i},$ for all $i\in I.$ Thus,  by the continuity of $m^x$ and $m^{x_0}$, we have that $m(g_i,x)\rightarrow m(g,x)$ and
	$m(g_i,x_0)\rightarrow m(g,x_0),$ which gives 
	\begin{equation*}\label{eqG1}
	(\forall x\in \X_{g^{-1}})\left(d(x,x_0)<\frac{1}{2^n}\Rightarrow d\left(m(g,x),m(g,x_0)\right)\leq\frac{1}{2^l}\right),
	\end{equation*}
	and we conclude that $g\in F_{n,l}.$ Now, we check that 
	\begin{equation}\label{eqG^0}
	G^{x_0}=\bigcap\limits_{l}\bigcup\limits_{n} F_{n,l}.
	\end{equation}
	It is clear that  $G^{x_0}\supseteq \bigcap\limits_{l}\bigcup\limits_{n} F_{n,l}.$  Conversely, take  $g\in G^{x_0}$  and $l\in \N$.  Since $m_g$ is continuous at $x_0$, for  $\varepsilon=\dfrac{1}{2^l}$, there exists  $\delta>0$ such that
	$$
	(\forall x\in \X_{g^{-1}})\left(d(x,x_0)<\delta\Rightarrow d\left(m(g,x),m(g,x_0)\right)\leq\frac{1}{2^l}\right),
	$$
	Let  $n\in \mathbb{N}$ with  $\dfrac{1}{2^n}<\delta$, then 
	$g\in F_{n,l}.$ Thus $G^{x_0}\subseteq\displaystyle\bigcap_l\bigcup_n F_{n,l},$ and we have \eqref{eqG^0}.
	
	Since $F_{l,n}$ is closed in $G^{x_0}$, for all $n,l;$  the set $D=\bigcup\limits_{l}\bigcup\limits_{n} (F_{n,l}\setminus {\rm int}(F_{n,l}))$ is meager. But  $G^{x_0}$ is a non empty open set, since $G$ is Baire, there is $g_1\in G^{x_0}\setminus D $. We shall check that $m$ is continuous at $(g_1,x_0).$ Indeed, let $\{(h_\alpha, y_\alpha)\}_{\alpha\in \Lambda}\subseteq G*\X$ be a net converging to $(g_1,x_0).$  Take  $\varepsilon >0$ and $l\in \mathbb{N}$ such that $\dfrac{1}{2^{l-1}}<\varepsilon.$ By 
	\eqref{eqG^0},  there exists $n\in \mathbb{N}$ such that  $g_1\in F_{n,l}.$ Since  $g_1\not\in D$  then $g_1\in \text{int}(F_{n,l})$. The fact that $h_\alpha\rightarrow g_1$, implies that there is $\alpha_1\in \Lambda$ such that  $h_\alpha\in \text{int}(F_{n,l}),$ for all $\alpha\geq \alpha_1.$ Also, since $y_\alpha\rightarrow x_0,$ there exists $\alpha_2\in \N$ such that  $d(y_\alpha,x_0)<\dfrac{1}{2^n},$ for all $\alpha\geq \alpha_2.$ Additionally, by the continuity of  $m^{x_0}$, we have  that  $m(h_\alpha,x_0)\rightarrow m(g_1,x_0).$ Thus, there exists $\alpha_3\in \Lambda$ such that $d(m(h_\alpha,x_0),m(g_1,x_0)) <\dfrac{1}{2^l},$ for all $\alpha\geq \alpha_3.$
	
	Let $\alpha\geq \max\left\{\alpha_1,\alpha_2,\alpha_3\right\}$. Then $h_\alpha\in F_{n,l}$ and $d(y_\alpha,x_0)<\dfrac{1}{2^n}$, 
	hence $$ d(m(h_\alpha,y_\alpha),m(g_1,x_0))\leq d(m(h_\alpha,y_\alpha),m(h_\alpha,x_0))+d(m(h_\alpha,x_0),m(g_1,x_0))<\frac{1}{2^{l-1}}<\varepsilon.$$
	Thus, $m$ is  continuous at $(g_1,x_0).$
	
	Now, since  $x_0\in \X_{g\m_{1}}\cap  \X_{g\m_{0}},$ we have, by Proposition \ref{fam}(ii), that  $g_1\cdot x_0 \in  \X_{g_{1}}\cap  \X_{g_1g\m_{0}},$  then $(g_0g\m_{1})\cdot (g_1\cdot x_0)$ is defined and  $(g_0g\m_{1})\cdot (g_1\cdot x_0)=g_0\cdot x_0,$ by (PA2). That is
	\begin{equation}
	\label{compdeu} m(g_0,x_0)=m(g_0g\m_{1},m(g_1, x_0)).
	\end{equation}
	
	Finally take a net $\{(h_j, y_j)\}_{j\in J}$ in $G*\X$ converging to $(g_0,x_0).$ Then 
	$g_1g\m_{0}h_j\to g_1\in G^{x_0},$ and since $ G^{x_0}$ is open we assume that the net $\{g_1g\m_{0}h_j\}_{j\in J}$ is contained in $G^{x_0},$ thus  $x_0\in \X_{(g_1g\m_{0}h_j)\m},$ for all $j\in J.$ But $y_j\to x_0,$ then there is $j_0\in J$ such that $y_j\in \X_{(g_1g\m_{0}h_j)\m},$ for all $j\geq j_0.$ 
	Now $y_j\in \X_{(g_1g\m_{0}h_j)\m}\cap \X_{h\m_{j}},$ then $h_j\cdot y_j\in\X_{(g_1g\m_{0})\m}$ (by Proposition \ref{fam}(ii)). By  (PA2) we get 
	$$
	(g_1g\m_{0}h_j)\cdot y_j=(g_1g\m_{0})\cdot (h_j\cdot y_j)\in \X_{g_1g\m_{0}}.
	$$ 
	Thus, by (PA1), 
	$$
	h_j\cdot y_j=(g_0g\m_{1})\cdot [(g_1g\m_{0}h_j)\cdot y_j].
	$$
	Notice that  $(g_1g_0\m h_j,y_j)\to(g_1,x_0)$ . Then by continuity of $m_{g_0g\m_1}$ and the  continuity of $m$ at  $(g_1,x_0)$  one gets
	$$
	m(h_j,y_j)=(g_0g\m_{1})\cdot [(g_1g\m_{0}h_j)\cdot y_j]\to(g_0g\m_{1})\cdot(g_1\cdot x_0)=m(g_0g\m_{1},m(g_1, x_0)),
	$$
	and  by  \eqref{compdeu} we obtain $m(h_j,y_j)\to m(g_0,x_0).$  
	\endproof

	\begin{cor} Let $G$ be a topological Hausdorff  group, and $a\colon G\times \X\to \X$ an action of $G$  on a metric space  $\X.$ If $G$ is Baire, then $a$ is continuous, if and only if, $a$ is separately continuous.
	\end{cor}

	\begin{cor} Let $G$ be a countable discrete group and $m$ a  topological partial action of $G$ on a metric space  $\X.$ Then $m$ is continuous, if and only if, $m$ is separately continuous.
	\end{cor}
	
%

\begin{exe}\label{mob}{\bf M{\"o}bius transfromations}  \cite[p. 175]{CHL} The group $G={\rm GL}(2,\mathbb{R})$ is Polish and acts partially on $\mathbb{R}$ by setting 

$$g\cdot x=\frac{ax+b}{cx+d},\,\,\,\,\,\,{\rm where}\,\,\,\,\,\, g= \left( \begin{array}{ccc}
a & b \\
c & d \\ \end{array} \right) \in G.$$
\end{exe}
Notice that for all $g\in G$ the set  $\mathbb{R}_g=\{x\in\mathbb{R}\mid cx+d\neq 0\}$ is open, and $m=\{m_g\colon \mathbb{R}_{g\m}\ni x\to g\cdot x\in \mathbb{R}_g\}_{g\in G}$ is a topological partial action.  For $x\in \mathbb{R} $ let $t_x=\left(\begin{array}{ccc}
1 & x \\
0 & 1 \\ \end{array}\right) ,$ then for $y\in \mathbb{R}$ one has that $t_{y-x}\cdot x=y.$ 
Since
$$G^0=\left\{\left(\begin{array}{ccc}
a & b \\
c & d \\ \end{array}\right) \in G: d\neq 0\right\}$$ is open,  and 

  $$G^x=G^{t_x\cdot 0}=G^0 t\m_x=G^0\left(\begin{array}{ccc}
1 &- x \\
0 & 1 \\ \end{array}\right),$$  then $G^x$ is open. Finally, since

$$m^x\colon G^x\ni\left(\begin{array}{ccc}
a & b \\
c & d \\ \end{array} \right) \mapsto\displaystyle \frac{ax+b}{cx+d}\in \mathbb{R}, $$ is continuous, then $m$ is continuous.

\medskip

\noindent\textbf{Acknowledgements}

\noindent The authors thank the referee for his/her suggestions and comments that helped to improve the quality of the paper. In particular, for suggesting the current version of Theorem \ref{sep}, that was initially established only for Polish groups. 

\end{document}